# Unique Structure on $G$-equivariant Manifolds

Reza Aghayan[*][†]

July 11, 2018

**Abstract**

We will have a deep look at the set of all $G$-equivariant maps from the factor Lie group $G$ to the under the action manifold $M$, both from "computational" and "observability" viewpoint. We will also be looking for the existence of "unique" structure on this set, in a way that the induced-action of Lie group $G$ be smooth on the mentioned set.

## 1 Introduction and Preliminaries

Moving back and forth between "computation" and "observation" has been always an impressive key in the history of science. Every thing that is being "observed" has some "computations" inside it and in its environment and every "computation" must be concluded with an applicable "observation". Here is, exactly, where mathematics, connecting these two ends, plays its fundamental role well and interprets the "observable facts" from the models derived from "computations".

---

[*]Kingston University, London, UK.
[†]E-mail: reza_raghayan@kingston.ac.uk

These kinds of, a manifold is an object which, locally, just looks like an open subset of some Euclidean space, but whose global topology can be quite different. More explicitly, an m-dimensional real manifold $M$ is a Hausdorff and second countable topological space which is covered by a collection of open subsets $\mathcal{W}_\alpha \subset M$ called coordinate charts, and one to one map $\chi_\alpha : \mathcal{W}_\alpha \to W_\alpha \subset \mathbb{R}^m$ onto a connected open subset, which serve to define local coordinates on $M$. The manifold is smooth if the transition functions $\chi_\alpha \circ \chi_\beta^{-1}$ are smooth where defined.

It is fitting to mention, any manifold $M$, as a topological space, is a"normal space", meaning that for each pair of disjoint closed subsets in $M$, there exist disjoint open subsets in $M$ containing the closed subsets, respectively.

As an another example, Lie groups, roughly speaking, which are smooth manifolds that are also "groups" in which multiplication and inversion are smooth maps.

Now we provide a rapid survey of some basic preliminaries based on the concept of Lie groups including some fundamental properties of them as group transformations and $G$-equivariant maps which will play a vital role in our studies in the following sections.

**Proposition 1.1.** *Let $G$ be Lie group and $g \in G$. Let $L_g : h \mapsto gh$, $R_g : h \mapsto hg$ and $C_g : h \mapsto ghg^{-1}$ denote, respectively, the associated left, right and conjugate multiplication maps of $G$ on itself. It is readily seen*

$$\begin{aligned} L_{gh} &= L_g L_h, & R_{gh} &= R_h R_g, \\ C_{g^{-1}} &= (C_g)^{-1}, & C_g &= L_g \circ R_{g^{-1}} = R_{g^{-1}} \circ L_g. \end{aligned} \tag{1}$$

and also $C_g$ is a homomorphism. In addition to, these maps are, clearly, diffeomorphisms on Lie group $G$.

Of particular importance is relation between these two objects, when Lie groups of transformations act on manifolds. A "(left-)transformation group" acting on a smooth manifold $M$ is determined by a Lie group $G$ and smooth map $\Psi : G \times M \to M$ denoted by $(g, x) \mapsto g \cdot x$, which satisfies

$$e \cdot x = x, \quad g \cdot (h \cdot x) = (gh) \cdot x \qquad \text{for all} \quad x \in M, \quad g \in G.$$

Manifold $M$ endowed with a specific action of Lie group $G$ is called a (left or right) "$G$-space".

Let manifold $M$ be a left $G$-space. Any smooth map defined by

$$\phi_g : M \to M; \quad x \mapsto g \cdot x \quad \text{for all} \quad x \in M,$$

provided that $g \cdot x$ is defined, is a diffeomorphism, or $(\phi_g)^{-1} = \phi_{g^{-1}}$. Also, any map

$$\psi_x : G \to M; \quad g \mapsto g \cdot x \quad \text{for all} \quad g \in G,$$

where $g \cdot x$ is defined, is smooth and its range, $O^x = G \cdot x = \{g \cdot x \mid g \in G\}$, is the orbit of $G$ through point $x$.

It is seen easily

$$\begin{aligned} i\,) & \quad \psi_{g \cdot x} = \psi_x \circ R_g, \quad \text{and} \\ ii\,) & \quad \psi_x \circ L_g = \phi_g \circ \psi_x. \end{aligned} \qquad (2)$$

**Definition 1.2.** *Suppose that manifolds $M$ and $N$ are left $G$-spaces. A smooth map $f : M \to N$ is said to be "$G$-equivariant", with respect to the given action $G$, if*

$$f(g \cdot x) = g \cdot f(x) \quad \text{for all} \quad g \in G.$$

**Proposition 1.3.** *Let Lie group $G$ act on smooth manifold $M$. Then set*

$$\{\psi_x : g \mapsto g \cdot x \mid x \in M\},$$

*is exactly the same set of all $G$-equivariant maps from Lie group $G$ to $M$, with respect to the left action $G$ on itself.*

**Proof.** Suppose $h \in G$, then

$$\psi_x(gh) = gh \cdot x = g \cdot \psi_x(h).$$

On the other hand, let $f : G \to M$ is a $G$-equivariant map. We claim $f = \psi_x$ where $x = f(e)$ :

$$\text{f}(h) = \text{f}(he) = h \cdot \text{f}(e) = h \cdot x = \psi_x(h) \quad \text{for all} \quad h \in G.$$

From now on, we suppose that Lie group $G$ acts globally on manifold $M$ and set of all $G$-equivariant maps from Lie group $G$ to the manifold $M$ will be **denoted** by $GM = \{\psi_x : g \mapsto g \cdot x \mid x \in M\}$.

In the next section, with some computations, we will find a strong and considerable system on the set $GM$. And then in Section 3, we will be investigating this set theoretically, and this helps us in carefully analyzing manifold topology and smooth structure on $GM$.

## 2  Facts in Computations

A vector field $\mathfrak{v}$ on Lie group $G$ is called "left- (right-) invariant" if $dL_g(\mathfrak{v}) = \mathfrak{v}$ (resp. $dR_g(\mathfrak{v}) = \mathfrak{v}$) for all $g \in G$ and denoted by $\mathfrak{v}_L$ (resp. $\mathfrak{v}_R$). Moreover, suppose Lie group $G$ acts on manifold $M$. A vector field v on manifold $M$ is called "$G$-invariant" if $d\phi_g(\mathfrak{v}|_x) = \mathfrak{v}|_{g \cdot x}$ for all $g \in G$ and all $x \in M$.

Note that any left- (right-) invariant vector field $\mathfrak{v}_L$ (resp. $\mathfrak{v}_R$) is uniquely determined by its value at the identity $e \in G$

$$\mathfrak{v}_L|_g = dL_g(\mathfrak{v}|_e),$$
$$(\text{resp.} \quad \mathfrak{v}_R|_g = dR_g(\mathfrak{v}|_e)).$$

Moreover,

$$\begin{aligned}\mathfrak{v}_L|_g &= \frac{d}{dt}|_{t=0} [\exp(t\mathfrak{v})g], \\ (\text{resp.} \quad \mathfrak{v}_R|_g &= \frac{d}{dt}|_{t=0} [g\exp(t\mathfrak{v})].\end{aligned} \quad (3)$$

where the set $\{\exp(t\mathfrak{v}) \mid t \in R\} \subset G$ refers to the "one-parameter" subgroup associated with a given tangent vector $\mathfrak{v}|_e \in T_eG$. Note that although the left- and right-invariant vector fields associated with a given tangent vector $\mathfrak{v}|_e \in T_eG$, are (usually) different and have different flows, nevertheless the associated one-parameter subgroups coincide. In other words

$$\exp(t\mathfrak{v}_R) = \exp(t\mathfrak{v}_L).$$

and denoted by $\exp(t\mathfrak{v})$.

Note further that, the identities (3) state that the left-invariant vector fields form the "infinitesimal generators" of the action of Lie group $G$ on $G$ by right multiplication and, vice versa, the "infinitesimal generators" of the action of $G$ on $G$ by left multiplication is the Lie algebra of right-invariant vector fields. [1],[2]

Suppose that $\mathfrak{v}_R$ is right-invariant vector field associated with $\mathfrak{v}|_e \in T_e G$. Then

$$\begin{aligned} \mathfrak{v}_L|_g &= dL_g(\mathfrak{v}|_e) = dL_g(dR_{g^{-1}}(\mathfrak{v}_R|_g)) \\ &= d(L_g \circ R_{g^{-1}})(\mathfrak{v}_R|_g) = d\mathrm{C}_g(\mathfrak{v}_R|_g). \end{aligned}$$

Therefore,

**Lemma 2.1.** *Suppose that $\mathfrak{v}_R$ and $\mathfrak{v}_L$ are right and left invariant vector fields associated to the tangent vector $\mathfrak{v}|_e \in T_e G$, then*

$$\begin{aligned} \mathfrak{v}_L|_g &= d\mathrm{C}_g(\mathfrak{v}_R|_g), \qquad \text{and} \\ \mathfrak{v}_R|_g &= d\mathrm{C}_{g^{-1}}(\mathfrak{v}_L|_g) = d\mathrm{C}_g^{-1}(\mathfrak{v}_L|_g). \end{aligned} \tag{4}$$

It is not also hard to check that

$$\begin{aligned} i\,) \quad & R_h \circ \mathrm{C}_g = \mathrm{C}_g \circ R_{g^{-1}hg}, \qquad \text{and} \\ ii\,) \quad & L_h \circ \mathrm{C}_g = \mathrm{C}_g \circ L_{g^{-1}hg}. \end{aligned}$$

Then

**Proposition 2.2.** *Let $\mathrm{C}_g : G \to G$, $g \in G$ be conjugate map, then the differential $d\mathrm{C}_g$ preserves the right and left invariance of vector fields.*

**Proof.** Let $\mathfrak{v}_R$ be right-invariant vector field on Lie group $G$. Then for any $h \in G$

$$\begin{aligned} dR_h(d\mathrm{C}_g(\mathfrak{v}_R)) &= d(R_h \circ \mathrm{C}_g)(\mathfrak{v}_R) \\ &= d(\mathrm{C}_g \circ R_{g^{-1}hg})(\mathfrak{v}_R) \\ &= d\mathrm{C}_g(dR_{g^{-1}hg}(\mathfrak{v}_R)) \\ &= d\mathrm{C}_g(\mathfrak{v}_R). \end{aligned}$$

Suppose Lie group $G$ acts on manifold $M$ and one-parameter subgroup $\{\exp(t\mathfrak{v}) \mid t \in R\} \subset G$ corresponds to the right-invariant vector field $\mathfrak{v}_R$ on Lie group $G$. Define

$$\mathrm{v}_R|_x = \frac{d}{dt}|_{t=0} [\exp(t\mathfrak{v}) \cdot x] \quad \text{for all} \quad x \in M. \tag{5}$$

In other words

$$\begin{aligned} \mathrm{v}_R|_x &= \frac{d}{dt}|_{t=0} [\exp(t\mathfrak{v}) \cdot x] \\ &= \frac{d}{dt}|_{t=0} [\psi_x(\exp(t\mathfrak{v}))] \\ &= d\psi_x(\frac{d}{dt}|_{t=0} [\exp(t\mathfrak{v})]), \end{aligned}$$

and hence

$$\mathrm{v}_R|_x = d\psi_x(\mathfrak{v}|_e), \tag{6}$$

defines a global vector field on manifold $M$. Since (1.1-i), we have

$$\begin{aligned} d\psi_x(\mathfrak{v}_R|_g) &= d\psi_x(dR_g(\mathfrak{v}|_e)) \\ &= d(\psi_x \circ R_g)(\mathfrak{v}|_e) \\ &= d\psi_{g \cdot x}(\mathfrak{v}|_e) = \mathrm{v}_R|_{g \cdot x}. \end{aligned}$$

In fact, the identity $\mathrm{v}_R = d\psi_x(\mathfrak{v}_R)$ together with next lemma illustrate that the differential $d\psi_x : TG \to TO^x \subset TM$, $x \in M$ gives an "alternative" definition for (5), but just when we consider it on the orbit $O^x$.

**Lemma 2.3.** *Suppose that Lie group $G$ acts on manifold $M$ and $\mathfrak{v}_R$ is a right-invariant vector field on $G$. Then the alternative definition*

$$\mathrm{v}_R|_{g \cdot x} = d\psi_x(\mathfrak{v}_R|_g), \quad g \in G, \quad x \in M, \tag{7}$$

*where $\psi_x : G \to M$ is a $G$-equivariant map, is independent of the choice of $x$. In other words, if $\psi_y : G \to M$ is another $G$-equivariant map corresponding to a point $y \in M$ so that $\mathrm{v}_R|_{g \cdot y} = d\psi_y(\mathfrak{v}_R|_g)$, then*

i) If $y \notin O^x$ then the differentials $d\psi_x$ and $d\psi_y$ have different ranges, respectively, $O^x$ and $O^y$, and there is nothing to check out.

ii) If $y \in O^x$ then
$$d\psi_x(\mathfrak{v}_R) = d\psi_y(\mathfrak{v}_R).$$

**Proof.** In case ii), it is equivalent to

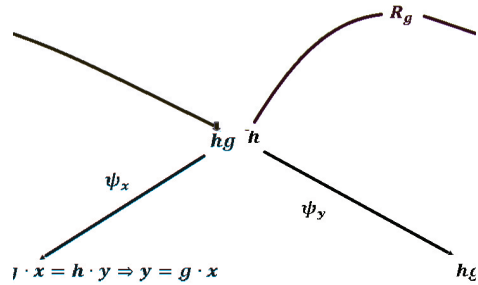

Then
$$\begin{aligned}
d\psi_y(\mathfrak{v}_R|_h) &= d\psi_{g \cdot x}(\mathfrak{v}_R|_h) \\
&= d\psi_x(dR_g(\mathfrak{v}_R|_h)) \\
&= d\psi_x(\mathfrak{v}_R|_{hg}).
\end{aligned}$$

Or,

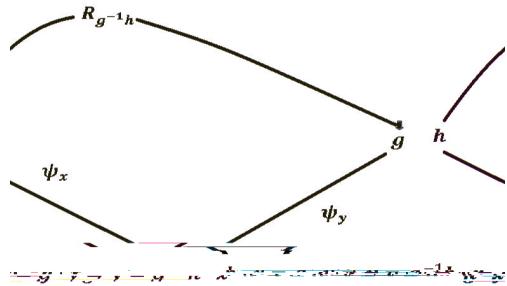

Then
$$\begin{aligned}
d\psi_y(\mathfrak{v}_R|_g) &= d\psi_{g^{-1}h \cdot x}(\mathfrak{v}_R|_g) \\
&= d\psi_x(dR_{g^{-1}h}(\mathfrak{v}_R|_g)) \\
&= d\psi_x(\mathfrak{v}_R|_h).
\end{aligned}$$

Therefore, the differential $d\psi_x : TG \to TO^x \subset TM$ defines, alternatively, the smooth vector fields $\mathrm{v}_R|_{O^x}$ on manifold $M$, given by (2.3), associated with a given right-invariant vector filed $\mathfrak{v}_R$ on Lie group $G$ and this definition is independent of the $G$-equivariant map $\psi_x$, $x \in M$. Note that the vector field $\mathrm{v}_R$ on manifold $M$ is **not** necessarily $G$-invariant.

A transformation group $G$ is said to be "free" if the only element of $G$ that fixes any element of $M$ is the identity, in other words, the isotropy group $G_x = \{e\}$ for all $x \in M$, equivalently, for $e \neq g \in M$, we have $g \cdot x \neq x$ for any $x \in M$. Besides, it is not hard to check that the action is free if and only if for each $x \in M$, the $G$-equivariant map $\psi_x : G \to M$ is injective.

Moreover, in general, the orbits of a Lie group of transformations are all submanifolds of the manifold $M$. Also, a group action of Lie group $G$ is called "regular" if its orbits are all submanifolds having the same dimension as submanifolds of $M$ and also for each point $x \in M$ there exists arbitrarily small neighbourhood $\mathcal{U}$ of $x$ with the property that each orbit of $G$ intersects $\mathcal{U}$ in a connected subset.

Note that if $G$ acts regularly on manifold $M$ then each orbit of $G$ is a regular (embeded) submanifold of $M$. The condition that each orbit be a regular submanifold is necessary, but not sufficient to guarantee regularity of the group action.

**Proposition 2.4** *An $r$-dimensional Lie group $G$ acts locally freely on a manifold $M$ if and only if its orbits have the same dimension $r$ as $G$ itself.*

In fact, the orbits of regular group actions have a simple local canonical form, generalizing the rectifying coordinates for a one-parameter group of transformations.

**Theorem 2.5** *Suppose Lie group $G$ acts regularly on an $m$-dimensional manifold $M$ with orbits of dimension $s$. Then, near every point of $x \in M$, there exist "rectifying local coordinates"*

$$x = (y, z) = (y^1, \ldots, y^s; z^1, \ldots, z^{m-s}),$$

*having the property that any orbit intersects the coordinate chart in at most one slice $\mathcal{N}_c = \{z^1 = c_1, \ldots, z^{m-s} = c_{m-s}\}$, for constants $c = (c_1, \ldots, c_{m-s})$.*

In fact, if action of a transformation group is regular then Theorem 2.5 demonstrates that the orbits form a "foliation" of the underlying manifold. [1]

**Remark.** It is worth mentioning that, in general, any complicated Lie group action can be "regularized" by lifting it to a suitable bundle sitting over the original manifold $M$, called "lifted action", such that the right and left regularizations of the transformation group $G$ define regular, free actions.

**Definition 2.6** *Let Lie group $G$ act on manifold $M$. The **left-lifted** action of the action of $G$ on $M$ is the "diagonal action" of $G$ on $\bar{M} = G \times M$ provided by the map*

$$\bar{L}_g(h, x) = \bar{L}(g, (h, x)) = (gh, g \cdot x), \quad g \in G, \quad (h, x) \in \bar{M}. \qquad (8)$$

*Also, the **right-lifted** action of $G$ is given by*

$$\bar{R}_g(h, x) = \bar{R}(g, (h, x)) = (hg^{-1}, g \cdot x), \quad g \in G, \quad (h, x) \in \bar{M}. \qquad (9)$$

These actions are called as the left or right **lifted action** of lie group $G$ since either projects back to the given action on $M$ via the $G$-equivariant projection $\pi_M : \bar{M} \to M$. In particular, for each point $(h, x) \in \bar{M}$, the diagram commutes, which maps $\psi_{(h,x)} \in G\bar{M}$ and $\psi_x \in GM$ are $G$-equivariant maps corresponding to the point $(h, x) \in \bar{M}$. The key, elementary result is that regularizing any group action eliminates all singularities and irregularities, e.g., lower dimensional orbits, non-embedded orbits, etc.

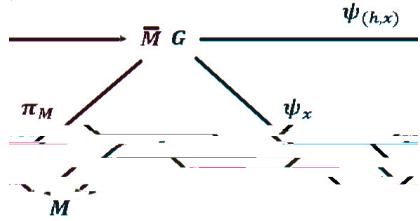

**Theorem 2.7.** *The right- and left-lifted actions of any transformation group $G$ define regular, free actions on the bundle $\bar{M} = G \times M$.*

Thus, lifting the action of $G$ on $M$ to the bundle $\bar{M}$ has the effect of completely eliminating any irregularities appearing in the original action.

In addition to, the orbits of Lie group $G$ in manifold $M$ are the projections of their lifted counterparts in $\bar{M}$ and all of the lifted orbits have the same dimension as $G$ itself. [4]

From now on, we **suppose** the given action $G$ on manifold $M$ is free, equivalently, the $G$-equivariant maps $\psi_x$ for all $x \in M$ are injective.

**Definition 2.8.** *Let $\mathfrak{v}$ be a (right or left) invariant vector field on Lie group $G$. Any (respectively, right or left) vector field $\mathfrak{w}$ on $G$ is called a g-**conjugate** of $\mathfrak{v}$ if*

$$\mathfrak{w}|_e = d\mathrm{C}_g(\mathfrak{v}|_e).$$

The unique g-conjugate vector field corresponding to each vector field $\mathfrak{v}$ will be denoted by $\mathfrak{v}^g$.

**Proposition 2.9.** *Let $\mathfrak{v}_R$ and $\mathfrak{v}_L$ be, respectively, right and left invariant vector fields on Lie group $G$ provided by a tangent vector $\mathfrak{v}|_e \in T_e G$. Also, let $\mathfrak{v}_R^g$ and $\mathfrak{v}_L^g$, respectively, be their corresponding g-conjugates provided by $\mathfrak{v}^g|_e = dC_g(\mathfrak{v}|_e)$. Then*

$$\begin{aligned} \mathfrak{v}_R^g|_g &= d\mathrm{C}_g(\mathfrak{v}_R|_g) = \mathfrak{v}_L|_g \qquad \text{and} \\ \mathfrak{v}_L^g|_g &= d\mathrm{C}_g(\mathfrak{v}_L|_g). \end{aligned} \qquad (10)$$

Suppose Lie group $G$ acts on manifold $M$. Now we draw inspiration from alternative definition (2.5) for transporting left-invariant vector fields on Lie group $G$ to manifold $M$.

**Definition 2.10.** *Suppose that Lie group $G$ acts on manifold $M$ and $\psi_x : G \to M$, $x \in M$ is a $G$-equivariant map. Let $\mathfrak{v}_L$ be a left-invariant vector field on Lie group $G$. Define*

$$\mathrm{v}^x = d\psi_x(\mathfrak{v}_L).$$

Since any $G$-equivariant map $\psi_x$, $x \in M$ is injective then the vector field $\mathrm{v}^x$ defined by Definition 2.10, is well-defined, but just on the orbit $O^x \subset M$. Also given point $x$ is called the "base point".

**Proposition 2.11.** *Any vector field $\mathrm{v}^x$, defined by Definition 2.10, is a $G$-invariant vector field on manifold $M$. In the other words*

$$\mathrm{v}^x|_{g\cdot x} = d\phi_g(\mathrm{v}^x|_x). \tag{11}$$

*Equivalently,*

$$\mathrm{v}^x|_{g\cdot x} = \frac{d}{dt}|_{t=0}\, [g\exp(t\mathfrak{v})\cdot x].$$

**Proof.** By Definition 2.10

$$\begin{aligned}
\mathrm{v}^x|_{g\cdot x} &= d\psi_x(\mathfrak{v}_L|_g) \\
&= d\psi_x(dL_g(\mathfrak{v}|_e)) \\
&= d(\psi_x \circ L_g)(\mathfrak{v}|_e).
\end{aligned}$$

By (1.1-ii), $\psi_x \circ L_g = \phi_g \circ \psi_x$, then

$$\begin{aligned}
\mathrm{v}^x|_{g\cdot x} &= d(\phi_g \circ \psi_x)(\mathfrak{v}|_e) \\
&= d\phi_g(d\psi_x(\mathfrak{v}|_e)) \\
&= d\phi_g(\mathrm{v}^x|_x).
\end{aligned}$$

Now we want to compute some results. As we mentioned earlier, that is enough we discuss only in the case of $y \in O^x$ as a "base point".

In the first step

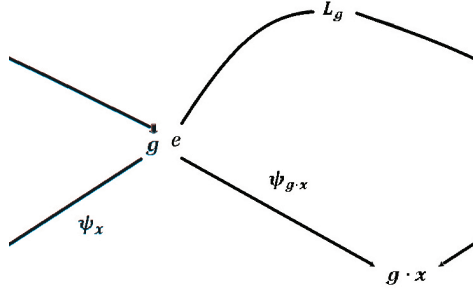

Then,

$$\begin{aligned}
d\psi_x(\mathfrak{v}_L|_g) &= d\psi_x(dL_g(\mathfrak{v}|_e)) \\
&= d(\psi_x \circ R_g)(dC_g(\mathfrak{v}|_e)) \\
&= d\psi_{g\cdot x}(\mathfrak{v}^g|_e).
\end{aligned}$$

And in the case of

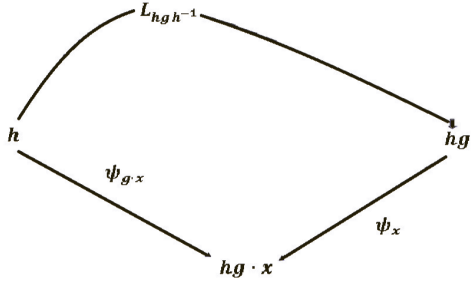

Then,
$$\begin{aligned} d\psi_x(\mathfrak{v}_L|_{hg}) &= d\psi_x(dL_{hgh^{-1}}(\mathfrak{v}_L|_h))) \\ &= d(\psi_x \circ R_g)[d(L_h \circ C_g \circ L_{h^{-1}})(\mathfrak{v}_L|_h))] \\ &= d\psi_{g \cdot x}(\mathfrak{v}_L^g|_h). \end{aligned}$$

Or
$$\begin{aligned} d\psi_{g \cdot x}(\mathfrak{v}_L|_h) &= d(\psi_x \circ R_g)(dL_h(\mathfrak{v}|_e)) \\ &= d\psi_x[(dL_{hg} \circ dC_{g^{-1}})(\mathfrak{v}|_e)] \\ &= d\psi_x[dL_{hg}(\mathfrak{v}^{g^{-1}}|_e)] = d\psi_x(\mathfrak{v}_L^{g^{-1}}|_{hg}). \end{aligned}$$

And in the general case of

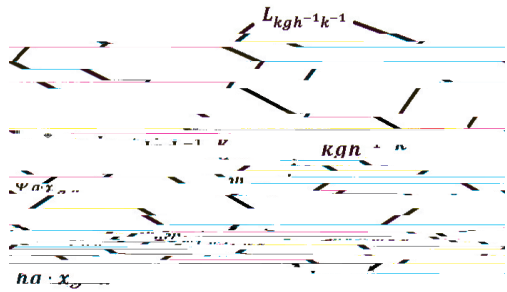

Then,

$$
\begin{aligned}
d\psi_{h\cdot x}(\mathfrak{v}_L|_{kgh^{-1}}) &= d(\psi_{g\cdot x} \circ R_{g^{-1}} \circ R_h)(dL_{kgh^{-1}k^{-1}}(\mathfrak{v}_L|_k)) \\
&= [d\psi_{g\cdot x} \circ d(R_{hg^{-1}} \circ L_k \circ L_{gh^{-1}} \circ L_{k^{-1}})](\mathfrak{v}_L|_k) \\
&= d\psi_{g\cdot x}[d(L_k \circ C_{gh^{-1}} \circ L_{k^{-1}})(\mathfrak{v}_L|_k)] \\
&= d\psi_{g\cdot x}(\mathfrak{v}_L^{gh^{-1}}|_k).
\end{aligned}
$$

**Theorem 2.12.** *Suppose Lie group $G$ acts on manifold $M$ and $\mathfrak{v}_L$ is a left-invariant vector field on $G$. Also suppose that $\mathrm{v}^x = d\psi_x(\mathfrak{v}_L)$ is the $G$-invariant vector field on manifold $M$, as in Definition 2.10, and $\psi_y : G \to M$ is another $G$-equivariant map corresponding to a point $y \in M$ and $\mathrm{v}^y = d\psi_y(\mathfrak{v}_L)$. Then*

i) *If $y \notin O^x$ then the differentials $d\psi_x$ and $d\psi_y$ have different ranges, respectively, $O^x$ and $O^y$, and there is nothing to check out.*

ii) *If $y \in O^x$ so that $y = g \cdot x$, then*

$$
\begin{aligned}
\mathrm{v}^x|_{hg\cdot x} &= (\mathrm{v}^g)^{g\cdot x}|_{hg\cdot x}, \qquad \text{or} \\
\mathrm{v}^{g\cdot x}|_{hg\cdot x} &= (\mathrm{v}^{g^{-1}})^x|_{hg\cdot x}.
\end{aligned}
$$

In fact any change in **base point** is equivalent to change to **conjugate vector fields**.

## 3  Facts in Observations

First of all, we provide a rapid survey of some basic propositions that will be needed in our studies in this section. [1],[5]

The "rank" of a map $f : M \to N$ at a point $x \in M$ is the rank of the Jacobian matrix $[\partial f_i / \partial x^j]$ of any local coordinate interpretation for the map f at the point $x \in M$. An "immersion" is a smooth function $f : M \to N$ with property that its differential, $df$, is injective at each point of manifold $M$, equivalently $\text{rank} f = \dim M$. A smooth "embedding" $f : M \to N$ is an injective immersion that is also a homeomorphism onto its image $f(M) \subset N$ in the subspace topology.

**Proposition 3.1.** *Let $M$ and $N$ be smooth manifolds. Let $f : M \to N$ be a smooth map of constant rank. If the map f is injective then it is an*

*immersion.*

**Proposition 3.2.** *Let* $f : M \to N$ *be a bijective smooth map of constant rank then it is a diffeomorphism.*

**Theorem 3.3.** *Suppose that $M$ and $N$ are smooth manifolds of the same dimension. Let $f : M \to N$ be a smooth immersion then $f$ is a local diffeomorphism. Moreover, a map $f : M \to N$ is a diffeomorphism if and only if it is a bijective local diffeomorphism.*

**Theorem 3.4.** *Suppose $M$ and $N$ are smooth manifolds and $G$ is a Lie group. Also let $f : M \to N$ be an $G$-equivariant map with respect to a smooth transitive action of $G$ on $M$ and any arbitrary smooth action $G$ on $N$. Then $f$ has constant rank. Moreover, its level sets are closed embeded submanifolds of $M$.*

It is fitting to mention Theorem 3.4 has considerable applications. As a straight consequence, every isotropy subgroup $G_x = \psi_x^{-1}\{x\}$ is a closed embeded submanifolds of Lie group $G$ and hence is a Lie subgroup.

As we mentioned earlier, each map $\psi_x : g \mapsto g \cdot x$ is $G$-equivariant with respect to the left multiplication on Lie group $G$ and the given action on manifold $M$. Also Lie group $G$ acts transitively on itself. Thus, by Theorem 3.4, any $G$-equivariant map $\psi_x : G \to M$ has constant rank and hence, by Proposition 3.1, is an immersion. Moreover, if the action is also regular then the orbits are all embeded submanifolds of $M$, thus the map $\psi_x : G \to M$ is an embedding. On the other hand, each map $\psi_x : G \to O^x$ is a bijective map from Lie group $G$ to the submanifold $O^x \in M$ (since $y = g \cdot x$ if and only if $\psi_x(g) = y$, thus $\psi_x : G \to O^x$ is surjective). Therefore, by Proposition 3.2, any map $\psi_x : G \to O^x$ is a diffeomorphism.

**Proposition 3.5.** *Suppose that Lie group $G$ acts freely on manifold $M$ and $\psi_x : G \to O^x \subset M$ belongs to $GM$. Then $\psi_x : G \to M$ is an immersion into manifold $M$ and a diffeomorphism from Lie group $G$ onto orbit $O^x$. Moreover, if the action of $G$ on $M$ is also regular then the map $\psi_x : G \to M$ is an embedding.*

Let $E$ and $M$ be smooth manifolds and $\pi : E \to M$ is a "vector bundle" of rank $k$. Let $\sigma^0 : M \to E$ given by

$$\sigma^0(x) = 0 \quad \text{for each} \quad x \in M,$$

denotes the zero-section of the vector bundle.

The zero-section of any smooth vector bundle is smooth. And also the image $\sigma(M)$ of a section is a smooth m-dimensional submanifold of $E$ that intersects each fiber $E_x = \pi^{-1}\{x\}$ in only one point.

If $F \subset E$ are vector bundles over manifold $N$ then the quotient bundle $E/F$ is a vector bundle with fibers $E_x/F_x$. If also $E$ is endowed with a fiberwise metric then we can identify $E/F$ with $E^\perp$.

In the case of $N \subset M$, take for $F \subset E$ to be $TN \subset TM|_N$ where $TM|_N$ is just the ambient tangent bundle restricted to points of $N$. In other words, the fiber of $E = TM|_N \to N$ at point $x$ is $T_xM$. Then the quotient bundle $E/F = (TM|_N)/TN$ is, by definition, the "normal bundle" to $N$.

Note that a Riemannian metric on manifold $M$ induces fiber metrics. Then we have a geometric picture of the "normal bundle" to $N$, denoted by $N^\perp$ as the set of vectors based at $N$, pointing "orthogonal" to $N$. In particular, rank of vector bundle $N^\perp$ is the codimension of $N$ in $M$.

For construction of normal bundle, first introduce a Riemannian metric on manifold $M$ so that we can identify $N_x^\perp$ with $(T_xN)^\perp$. Now based in the theory of geodesics, it can be asserted that through each point $(x, v) \in TM$ there is a unique curve $\gamma$, called a "geodesic", which playing the role of a straight line relative to the metric given, and satisfying the initial condition $\gamma(0) = x$, $\ddot{\gamma}(0) = v$. Denote

$$\gamma(t) = \exp(tv)x.$$

We can restrict the exponential map to the normal bundle

$$\exp|_{N^\perp} : N^\perp \to M,$$

by sending $(x, w) \in N^\perp$ to $\exp(w)x$. Under this map a line $tw$ in the fiber $N_x^\perp = (T_xN)^\perp$ gets sent to the unique geodesic tangent to $w$ and hence "orthogonal" to $N$ and passing through $x \in N$.

**Theorem 3.6.** *(Tubular Neighbourhoods Theorem) Let $M$ be a manifold and $N \subset M$ be an embeded submanifold. Then the map*

$$\exp|_{N^\perp} : N^\perp \to M \quad \text{given by} \quad (x, w) \mapsto \exp(w)x,$$

maps a neighbourhood of the zero-section in $N^\perp$ diffeomorphically onto a neighbourhood of an embeded submanifold $N$ in manifold $M$.

The Tubular Neighbourhood Theorem, in fact, asserts that sufficiently small neighbourhoods of an embeded submanifold $N$ in manifold $M$ are diffeomorphic to the normal bundle, denoted by $N^\perp$, in such a way that under this diffeomorphism submanifold $N$ gets mapped to the zero-section. In particular, by using a tubular neighbourhood, we can reduce questions in analysis near the submanifold to analysis in the normal bundle which is linear in the fiber.

By the Tubular Neighbourhoods Theorem, Theorem 2.5 and Proposition 3.5 we have the following fundamental theorem.

**Theorem 3.7.** *Let Lie group $G$ act freely and regularly on manifold $M$. Then we can introduce an adapted chart with the action $G$ on $M$, called **flat local coordinates***

$$x = (h, z) = (h^1, \ldots, h^r; z^1, \ldots, z^{m-r}), \quad h \in G, \quad z \in Z, \qquad (12)$$

*that locally identify $M$ with a subset of the Cartesian product $G \times Z$ with $Z \cong \mathbb{R}^{m-r}$ and such that the action of $G$ reduces to the trivial left action $g \cdot x = (g \cdot h, z)$.*

**Corollary 3.8.** *Suppose that Lie group $G$ acts freely and regularly on manifold $M$. Then each orbit of the action is a **closed** submanifold in $M$.*

**Definition 3.9.** *Suppose Lie group $G$ acts regularly on m-dimensional manifold $M$ with s-dimensional orbits. A (local) **cross-section** in $M$ is a (m-s)-dimensional submanifold $\mathcal{K} \subset M$ such that $\mathcal{K}$ intersects each orbit "transversally". The cross-section is regular if $\mathcal{K}$ intersects each orbit at most once.*

**Remark.** In flat local coordinates, a general cross-section is given by the graph $\mathcal{K} = \{(\varrho(z), z)\}$ of a smooth map $\varrho : Z \to G$. When we use flat local coordinates, we shall always assume, without loss of generality, that the identity cross-section $\{e\} \times Z$, i.e., when $\varrho(z) \equiv e$, belongs to the flat coordinate chart.

**Corollary 3.10.** *Let Lie group $G$ act freely and regularly on manifold $M$. Then for each $x \in M$ there is local regular cross-section $\mathcal{K}$ passing through point $x$ and intersects the orbit $O^x$ transversally.*

Suppose Lie group $G$ acts on a manifold $M$. The set of all orbits of $G$ in $M$, denoted by $M/G$, with the quotient topology should be called the **cross-section space** of the action. It is of great importance to determine conditions under which an orbit space is a smooth manifold.

**Lemma 3.11.** *Suppose that $M, Q$ and $P$ are smooth manifolds and map $\pi: M \to Q$ is a surjective submersion. Let $f: Q \to P$ is any map. Then $f$ is smooth if and only if $f \circ \pi$ is smooth.*

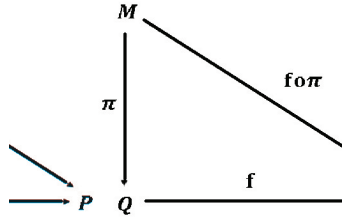

**Theorem 3.12.** *Let Lie group $G$ act smoothly, freely and regularly on manifold $M$. Then the cross-section space $M/G$ is a topological manifold of dimension equal to $dimM\text{-}dimG$, and has a unique smooth structure with the property that the quotient map $\pi: M \to M/G$ is a smooth submersion.*

**Proof.**
**Step 1.** *$M/G$ is a topological manifold.*

For any open set $\mathcal{U} \subset M$, we have

$$\pi^{-1}(\pi(\mathcal{U})) = \bigcup_{g \in G} \phi_g(\mathcal{U}).$$

Since $\phi_g: M \to M$ is a diffeomorphism, each such set is open and hence $\pi^{-1}(\pi(\mathcal{U}))$ is open in manifold $M$. Because $\pi: M \to M/G$ is a quotient map, this implies that $\pi(\mathcal{U})$ is open set in $M/G$, and therefore $\pi$ is an open map. If set $\{\mathcal{W}_\alpha\}$ is a countable basis for the topology of $M$ then set $\{\pi(\mathcal{W}_\alpha)\}$ is a countable collection of open subsets of $M/G$ and is a basis for the topology of $M/G$. Thus the cross-section space $M/G$ is second countable. On the other hand, since $\pi : M \to M/G$ is an open map then by Corollary 3.8, $M/G$ is a Hausdorff space.

In accordance with Proposition 3.5, since Lie group $G$ acts freely and regularly on manifold $M$ then the map $\psi_x : G \to M$ is an embedding and by Theorem 3.7, for each $x \in M$ there exists a flat local coordinates centered at $x$.

Now suppose that $q = \pi(x)$ is an arbitrary point of the space $M/G$. And also let $(\mathcal{U}, \chi)$ denote a flat local coordinates centered at point $x$, with $\chi(\mathcal{U}) = U_1 \times U_2 \subset R^r \times R^{m-r}$. Since $\pi : M \to M/G$ is an open map then subset $\mathcal{V} = \pi(\mathcal{U})$ is an open subset of $M/G$. Consider the coordinate functions of map $\chi$ in the form of $(y^1, \ldots, y^r; z^1, \ldots, z^{m-r})$, as in Theorem 3.7 and let $Z_0$ be the slice given by $\{y^1 = \ldots = y^r = 0\}$. It is not hard to see that $\pi : Z_0 \to \mathcal{V}$ is bijective by the properties of a flat local coordinates. Moreover, if $\tilde{\mathcal{V}}$ is an open subset of $Z_0$ then

$$\pi(\mathcal{W}) = \pi(\{(y, z) \mid (0, z) \in \mathcal{W}\})$$

is open in $M/G$ and thus $\pi|_{Z_0}$ is a homeomorphism. Suppose $\sigma : \mathcal{V} :\to Z_0$ to be the inverse of $(\pi|_{Z_0})^{-1}$ which is a local section of $\pi$.

Now, define a map $\theta = \pi_2 \circ \chi \circ \sigma : \mathcal{V} \to U_2$, where $\pi_2 : U_1 \times U_2 \to U_2$ is natural projection by sending the equivalence class of a point $(y, z)$ to $z$. This is well defined by the definition of a flat local coordinates. Because $\sigma$ and $\pi_2 \circ \chi$ are homeomorphisms, therefore $\theta$ is a homeomorphism. This completes the proof that $M/G$ is a topological manifold of dimension $m - r$.

**Step 2.** *$M/G$ is a smooth manifold.*

Finally, we need to show that $M/G$ has a smooth structure such that $\pi$ is a submersion. We will use the atlas consisting of all charts $(\mathcal{V}, \theta)$ as constructed in the preceding paragraph. With respect to any such chart for $M/G$ and the corresponding flat local coordinates for manifold $M$, the map $\pi : M \to M/G$ has the coordinate representation $\pi(y, z) = z$, which is certainly a submersion. Thus we need only show that any two such charts for $M/G$ are smoothly compatible.

Suppose $(\mathcal{U}, \chi)$ and $(\tilde{\mathcal{U}}, \tilde{\chi})$ are two flat local coordinates for manifold $M$ and let $(\mathcal{V}, \theta)$ and $(\tilde{\mathcal{V}}, \tilde{\theta})$ be the corresponding charts for $M/G$. Since $\phi_g$ is a diffeomorphism, we shall not lose any generality by restricting our attention to the case in which the two flat local coordinates are both centered at the same point $x \in M$. Writing the flat local coordinates as $(y, z)$ and $(\tilde{y}, \tilde{z})$, the fact that the coordinates are adapted to the action of $G$ means that two points with the same $z$-coordinate are in the same orbit, and therefore also have the same $\tilde{z}$-coordinate. This means that the transition map between these coordinates can be written $(\tilde{y}, \tilde{z}) = (\xi(y, z), \zeta(z))$, where $\xi$ and $\zeta$ are smooth functions defined on some neighbourhood of the origin. The transition map $\tilde{\theta} \circ \theta^{-1}$ is just $\tilde{z} = \zeta(z)$, which is clearly smooth.

The uniqueness of the smooth structure is an immediate result of the Lemma 3.11.

Suppose r-parameter Lie group $G$ acts on m-dimensional manifold $M$ and consider the set of all $G$-equivariant maps $GM = \{\psi_x : G \to M \mid x \in M\}$. Define action

$$G \times GM \to GM \quad \text{given by} \quad (g, \psi_x) \mapsto g \cdot \psi_x = \psi_{g \cdot x}.$$

From now on, we **call** the action of Lie group $G$ on manifold $M$ as "original-action" and the action of $G$ on the set $GM$ as "induced-action".

**Proposition 3.13.**
a) *The isotropy subgroups for original-action and induced-action are the same, meaning that, $G_{\psi_x} = G_x$. Moreover, the induced-action is transitive (free) if and only if the original-action is transitive (free).*

b) *Define map*

$$\psi : M \to GM \quad \text{given by} \quad x \mapsto \psi_x, \qquad (13)$$

*from manifold $M$ to the set of all $G$-equivariant maps $GM = \{\psi_x : G \to M\}$. This is a bijection and $G$-equivariant and also $O^{\psi_x} = \psi|_{O^x}$.*

**Proof.**
$$\begin{aligned} G_{\psi_x} &= \{g \in G \mid g \cdot \psi_x = \psi_x\} \\ &= \{g \in G \mid \psi_{g \cdot x} = \psi_x\} \\ &= \{g \in G \mid g \cdot x = x\} = G_x. \end{aligned}$$

And also

$$\begin{aligned} O^{\psi_x} &= G \cdot \psi_x = \psi_{G \cdot x} \\ &= \{\psi_y \mid y \in O^x\} = \psi|_{O^x}. \end{aligned}$$

Now the important question is, "Is there a unique manifold topology and smooth structure on the set of all $G$-equivariant maps $GM$ such that the induced-action is smooth?".

Clearly, by transforming the manifold topology and smooth structure of manifold $M$ to the set $GM$, so that the map $\psi : M \to GM$, defined in (3.2), be a diffeomorphism between manifolds, the set $GM$ is an m-dimensional manifold and diffeomorphic to the same manifold $M$ and the induced-action is smooth too. We call this structure on $GM$ as "induced-structure". Is this structure unique?

**Step I.** *Let the action Lie group $G$ on manifold $M$ be transitive.*

First we provide a rapid survey of some basic theorems which we will use in the following.

**Theorem 3.14.** *Suppose $G$ is a Lie group and let $H$ be a closed Lie subgroup of $G$. Then the left coset space $G/H$ is a smooth manifold, and the quotient map $\pi : G \to G/H$ is a smooth submersion.*

**Theorem 3.15.** *Let Lie group $G$ act transitively on manifold $M$ and $G_x$ denotes the isotropy Lie subgroup of point $x \in M$. Then the induced map $\bar{\psi}_x : G/G_x \to M$, given by $gG_x \mapsto g \cdot x$, is an $G$-equivariant diffeomorphism.*

**Theorem 3.16.** *Let $X$ be a set and Lie group $G$ act transitively on $X$ such that the isotropy group of a point $x \in X$ is a closed Lie subgroup of $G$. Then $X$ has a unique manifold topology and smooth structure such that the given action is smooth.*

By Theorem 3.14 the quotient space $G/G_x$ is a smooth manifold. What is important to observe is, the unique topology and smooth structure on the set $X$ is defined such that $\bar{\psi}_x : G/G_x \to X$ to be a diffeomorphism, as in the Theorem 3.15. Moreover, the given action of $G$ on $X$ is smooth because it can be written $(g, z) \mapsto \bar{\psi}_x(g \cdot \bar{\psi}_x^{-1}(z))$.

In other words, there is a unique manifold with smooth structure, up to diffeomorphisms, such that the action a given Lie group $G$ on it is smooth and transitive. If $H \subset G$ is the isotropy group of any arbitrary point of the manifold, this unique smooth manifold is usually denoted by $G/H$.

**Corollary 3.17.** *Let $G$ be a Lie group. Then there is a unique manifold topology and smooth structure on Lie group $G$ such that the left action $G$ on $G$ by $L_g$ is smooth. Moreover, for the left action of $G$ on $G$*

$$GG = \{L_g : G \to G \mid L_g(h) = g \cdot h\}$$

Note that the same corollary is true for the right action of $G$ on $G$.

**Corollary 3.18.** *suppose Lie group $G$ acts transitively on manifold $M$. Then there is a unique manifold topology and smooth structure on the set of all $G$-equivariant maps $GM$ such that the induced-action is smooth. It is equivalent to*

$$GM \cong M \cong G/H,$$

*where $H$ is isotropy group of an arbitrary point of the manifold $M$. And also since the diagram commutes, the map $\psi : M \to GM$ provided by $x \mapsto \psi_x$ is*

$$\begin{array}{ccc}
 & G/G_x & \\
\bar{\psi}_x \nearrow & & \nwarrow \bar{\psi}_{\psi_x} \\
M & \xrightarrow{\psi} & GM
\end{array}$$

*diffeomorphism. In particular,*

$$GG \cong G.$$

Now consider the restriction of original-action $\Psi : G \times M \to M$ and induced-action $\bar{\Psi} : G \times GM \to GM$ to any couple corresponding orbits, for instance $O^x \subset M$ and corresponding orbit $O^{\psi_x} \subset GM$. By Corollary 3.18, there is a unique manifold topology and smooth structure on the subset $O^{\psi_x}$ such that the restriction $\bar{\Psi}|O^{\psi_x} : G \times O^{\psi_x} \to O^{\psi_x}$ is smooth and the map $\psi|_{O^x} : O^x \to O^{\psi_x}$ is diffeomorphism.

**Lemma 3.19.** *Suppose Lie group $G$ acts smoothly on manifold $M$ and $\psi : x \mapsto \psi_x$ denotes the $G$-equivariant map from $G$ to the set $GM$. For any couple corresponding orbits, for instance $O^x \subset M$ and $O^{\psi_x} \subset GM$ there exists a unique manifold topology and smooth structure on $O^{\psi_x}$ such that the restriction of induced-action $\bar{\Psi}|O^{\psi_x} : G \times O^{\psi_x} \to O^{\psi_x}$ is smooth. Moreover, the restriction*

$$\psi|_{O^x} : O^x \to O^{\psi_x},$$

*is diffeomorphism. In particular*

$$O^x \cong O^{\psi_x} \cong G/G_x,$$

*where $G_x$ is the isotropy subgroup of the point $x \in M$.*

**Remark.** Suppose that there exists a manifold topology and smooth structure on the set $GM$ such that the induced-action of Lie group $G$ on manifold $M$ is smooth. It is obvious that the induced topology and smooth structure on each orbit $O^{\psi_x} \subset GM$ is the same as mentioned in Lemma 3.19.

**Lemma 3.20.** *Let $M, N$ and $P$ are smooth manifolds. Let $\mathrm{f} : M \to N$ be a surjective map from $M$ onto $N$. If $\mathrm{f} \circ \mathrm{g} : P \to N$ is a continues map then $\mathrm{g} : P \to M$ is also continues.*

**Proposition 3.21.** *Suppose that $\Psi : G \times M \to M$ denotes the smooth original-action of Lie group $G$ on manifold $M$ and $\bar{\Psi} : G \times M \to M$ denotes the induced-action $G$ on the set of all $G$-equivariant maps $GM$. Also, let there exist a manifold topology and smooth structure on the set $GM$. In this case,*
  i) *If the map $\psi : x \mapsto g \cdot x$ is smooth then the induced-action $G$ on $GM$ is smooth. In fact,*

$$\bar{\Psi} = \psi \circ \Psi.$$

  ii) *If the induced-action $G$ on $GM$ is smooth then the map $\psi : G \to GM$ is at least continues. Moreover, since $\psi^{-1} \circ \bar{\Psi} = \Psi$, then $\psi : G \to GM$ is at least homeomorphism.*

**Proof.**

$$\bar{\Psi}(g, \psi_x) = g \cdot \psi_x = \psi_{g \cdot x} = \psi(g \cdot x) = (\psi \circ \Psi)(g, x).$$

**Theorem 3.22.** *Suppose Lie group $G$ acts smoothly and regularly on manifold $M$. Let there exist a manifold topology and smooth structure on the set $GM$ such that the induced-action $\bar{\Psi} : G \times M \to M$ is smooth. Then the induced-action $G$ on manifold $GM$ is regular.*

**Proof.** On the one hand, because the original-action $G \times M \to M$ is regular so, by definition, its orbits are all submanifolds having the same dimension as submanifolds of manifold $M$ and hence by Lemma 3.19, all orbits of the induced-action $G \times GM \to GM$ are submanifolds having the same dimension as submanifolds of manifold $GM$.

On the other hand, for regularity of the induced-action we must also show that for each point $\psi_y \in GM$ there exists small neighbourhood $\bar{\mathcal{U}}$ of $\psi_y$ with the property that each orbit of $G$ in $GM$ intersects $\bar{\mathcal{U}}$ in a connected subset.

Suppose that for a point $\psi_y \in GM$ and any small neighbourhood $\bar{\mathcal{U}}$ of $\psi_y$ there is an orbit, namely $O^{\psi_x} \subset GM$ where $\psi_x \in \bar{\mathcal{U}}$, such that intersects $\bar{\mathcal{U}}$ in at least two connected components.

By Proposition 3.20-ii), since $\psi : G \to GM$ is a homeomorphism so $\psi^{-1}(\bar{\mathcal{U}})$ is a connected open subset of the corresponding point $x \in M$. On the other hand, because the map $\psi$ is continues (the image of a connected subset under a continues map is connected), its restriction of to each orbit is diffeomorphism and also the original-action is regular, these all leads to a contradiction.

**Step II.** *Let the original-action $G \times M \to M$ be smooth, free and regular then there exists a unique manifold topology and smooth structure on the set $GM$ such that the induced-action $G \times GM \to GM$ is smooth. In particular,*

$$GM \cong M.$$

Since the action of Lie group $G$ on manifold $M$ is smooth, free and regular and the induced-action of $G$ on manifold $GM$ is smooth, then by Proposition 3.13-a) and Lemma 3.21, the induced-action is free and regular.

Therefore, by Theorem 3.7, there exist local flat coordinates

$$(h, z) = (h^1, \ldots, h^r; z^1, \ldots, z^{m-r}), \quad h \in G, \quad z \in Z,$$

centered at each arbitrary point $x_0 \in M$ and

$$(\bar{h}, \bar{z}) = (\bar{h}^1, \ldots, \bar{h}^r; \bar{z}^1, \ldots, \bar{z}^{m-r}), \quad \bar{h} \in G, \quad \bar{z} \in Z,$$

centered at the corresponding point $\psi(x_0) = \psi_{x_0} \in GM$, that locally identify the manifolds with subsets of the Cartesian product $G \times Z$, with $Z \cong \mathbb{R}^{m-r}$, and the actions of Lie group $G$ on manifolds $M$ and $GM$, reduce to the trivial left action

$$g \cdot (h, z) \;=\; (g \cdot h, z),$$
$$g \cdot (\bar{h}, \bar{z}) \;=\; (g \cdot \bar{h}, \bar{z}), \qquad \text{respectively.}$$

By Proposition 3.20-ii), the map $\psi : G \to GM$ is homeomorphism thus we can consider small neighbourhoods of any coupled corresponding points $x_0 \in M$ and $\psi(x_0) = \psi_{x_0} \in GM$ with mentioned flat local coordinates. In these small neighbourhoods,

$$(\bar{h}, \bar{z}) = \psi(h, z) = \psi_{(h,z)} \in GM.$$

Now we represent the $G$-equivariant map $\psi : G \to GM$ in these flat local coordinates. Thus

$$\psi : G \times Z \to G \times Z; \quad (h, z) \mapsto \psi_{(h,z)}.$$

Denote $\psi = (\psi_1, \psi_2)$. We will be looking for, in fact, a representation for component maps $\psi_1 : G \times Z \to G$ and $\psi_2 : G \times Z \to Z$.

**Definition 3.23.** *Let Lie group $G$ act on a manifold $M$. A function $f : M \to N$, where $N$ is a manifold, is called an "Invariant function" if $f(g \cdot x) = f(x)$, where defined.*

**Lemma 3.24.** *Let $\psi = (\psi_1, \psi_2) : (h, z) \mapsto \psi_{(h,z)}$ be the same as mentioned above. Then $\psi_1 : G \times Z \to G$ is a $G$-equivariant map and $\psi_2 : G \times Z \to Z$ is an Invariant function.*

**Proof.** Since $\psi : G \times Z \to G \times Z$ is a $G$-equivariant map, so

$$\psi(g \cdot (h, z)) = g \cdot \psi(h, z).$$

In the other words

$$\begin{aligned}(\psi_1(g \cdot (h, z)), \psi_2(g \cdot (h, z))) &= g \cdot (\psi_1(h, z), \psi_2(h, z)) \\ &= (g \cdot (\psi_1(h, z)), \psi_2(h, z)).\end{aligned}$$

So

$$\begin{aligned}\psi_1(g \cdot (h, z)) &= g \cdot (\psi_1(h, z)) \quad \text{and} \\ \psi_2(g \cdot (h, z)) &= \psi_2(h, z) \qquad \text{for all} \quad (h, z) \in G \times Z.\end{aligned} \qquad (14)$$

And hence

$$\psi_2(g, z) = \psi_2(g \cdot (e, z)) = \psi_2(e, z) \quad \text{for all} \quad g \in G. \tag{15}$$

So the definition of $\psi_2 : G \times Z \to Z$ is independent of group elements, and hence we can replace it by a map $\bar{\psi}_2 : Z \to Z$ provided by

$$\bar{\psi}_2(z) = \psi_2(e, z) = \psi_2(g, z) \quad \text{for all} \quad g \in G.$$

It is not hard to show that $\bar{\psi}_2$ is one to one, onto and the law of the inverse map of $\psi$, in given local flat coordinates is in the form of

$$\psi^{-1} = (*, \bar{\psi}_2^{-1}).$$

In other words, the map $\bar{\psi}_2 : Z \to Z$ is a homeomorphism.

**Lemma 3.25.** *Suppose $\psi = (\psi_1, \psi_2) : (h, z) \mapsto \psi_{(h,z)}$ is the same as mentioned in Lemma 3.24. Then $\psi_1 : G \times Z \to G$ is a G-equivariant map and $\psi_2 : Z \to Z$ is a homeomorphism. Moreover,*

$$\psi|_{G \times \{0\}} = \mathbb{I}d.$$

These "computations" as well as some "intuitive ideas" which I will present those in Section 4, give me enough motivation to state the following conjecture.

**Conjecture** *Suppose Lie group $G$ acts smoothly, freely and regularly on manifold $M$ and $GM = \{\psi_x : x \mapsto g \cdot x \mid x \in M\}$ denotes all G-equivariant maps from $G$ to $M$. Also suppose $\psi : G \to GM$ is the G-equivariant map provided by $x \mapsto \psi_x$. Then $\psi$ is a **diffeomorphism**.*

Therefore, in the case Step II,

$$GM \cong M.$$

**Step III.** *Let Lie group $G$ act smoothly on manifold $M$. Then there is a unique manifold topology and smooth structure on the set $GM$ such that the induced-action $G$ on manifold $GM$ is smooth and also the map $\psi : M \to GM$ is smooth.*

By Definition 2.6, there is left-lifted action of the action of Lie group $G$ on manifold $\bar{M} = G \times M$ given by the map

$$\bar{L}_g(h,x) = \bar{L}(g,(h,x)) = (g \cdot h, g \cdot x), \quad g \in G, \quad (h,x) \in \bar{M}.$$

and respectively, there is left-lifted action of the action of Lie group $G$ on $G\bar{M} = G \times GM$ provided by the map

$$\bar{L}_g(h, \psi_x) = \bar{L}(g,(h,\psi_x)) = (g \cdot h, \psi_{g \cdot x}), \quad g \in G, \quad (h,\psi_x) \in G\bar{M}.$$

such that these action is free and regular.

**Proposition 3.25.** *Let Lie group $G$ act smoothly on manifold $M$ and $G \times \bar{M} \to \bar{M}$ refer to the corresponding left-lifted action, where $\bar{M} = G \times M$ denotes the lifted manifold. Also suppose $GM$ is the set of all $G$-equivariant maps from $G$ to $M$, then*

$$G\bar{M} = \{(R_h, \psi_x) : G \to \bar{M} \mid h \in G,\ x \in M\}.$$

Moreover, the following diagram commutes.

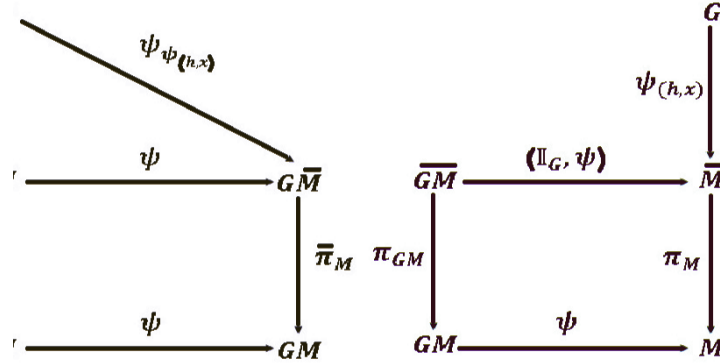

**Proof.** Let $\psi_{(h,x)} \in G\bar{M}$, then

$$\begin{aligned}\psi_{(h,x)}(g) &= g \cdot (h, x) \\ &= (g \cdot h, g \cdot x) \\ &= (R_h(g), \psi_x(g)).\end{aligned}$$

And hence

$$G\bar{M} = \{\psi_{(h,x)} : G \to \bar{M} \mid (h,x) \in \bar{M}\}$$
$$= \{(R_h, \psi_x) : G \to \bar{M} \mid h \in G,\ x \in M\}.$$

Therefore, by step II

$$G\bar{M} \cong \bar{M}.$$

It is equivalent to

$$\{(R_h, \psi_x) \mid h \in G,\ x \in M\} \cong G \times M,$$

and hence

$$GM = \{\psi_x : g \mapsto g \cdot x \mid x \in M\} \cong M,$$

# 4 Appendix

Unfortunately, since considered work is not in my recent major studies, I have to leave some computations and observations half done. The most important of them is finding a way to "Moving frame",[3],[4]. As far as some computations indicates (I have not present them here), the degree of non-smoothness in the above discussion is as difficult as the lack of $G$-equivariant function

$$\rho : M \to G.$$

Studying Lie algebras could be desired too.

Among them is the same which me to state the conjecture. As I noted above, in this part I will present my "intuitive ideas" which led me to the conjecture.

Since the map $\psi$ is transverse to orbits and its restriction to each orbit in manifold $M$ is diffeomorphism to a unique orbit in manifold $GM$ and also $\psi_2 : Z \to Z$ is a homomorphism, then we can replace the corresponding flat local coordinates at $\psi(x) = \psi_{x_0} \in GM$ by a new flat coordinates in the form of $G \times \psi_2(Z) \cong G \times Z$ at this point.

In fact, we choose the new flat local coordinates at point $\psi_{x_0}$ such that the coordinates of any corresponding point $(e, \psi_2(z)) \in GM$, in the new flat

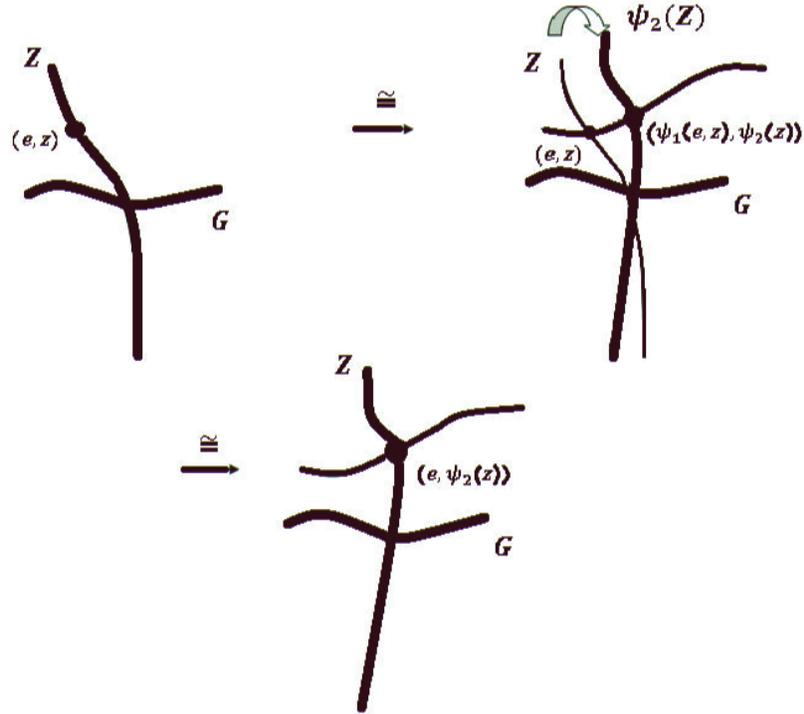

local coordinates, is the same as the coordinates of the point $(e, z) \in M$. Thus

$$\begin{aligned} \psi(g, z) &= (\psi_1(g, z), \psi_2(z)) \\ &= (g \cdot \psi_1(e, z), \psi_2(z)) \\ &\cong (g, z). \end{aligned}$$

Therefore, $\psi : x \mapsto \psi_x$ is a local diffeomorphism from manifold $M$ to the manifold $GM$. Moreover, since the $\psi$ is bijective, by Proposition 3.3, it is a diffeomorphism.